\documentclass[11pt]{amsart}

\newtheorem{theorem}{Theorem}[section]
\newtheorem{lemma}[theorem]{Lemma}

\theoremstyle{definition}
\newtheorem{definition}[theorem]{Definition}
\newtheorem{question}[theorem]{Question}
\newtheorem{example}[theorem]{Example}

\newtheorem{corollary}[theorem]{Corollary}
\newtheorem{remark}[theorem]{Remark}
\theoremstyle{remark}

\newcommand{\be}{\begin{equation}}
\newcommand{\ee}{\end{equation}}

\numberwithin{equation}{section}

\begin{document}

\title{Chern numbers and the indices of some elliptic differential operators}

\author{Ping Li}
\address{Department of Mathematics, Tongji University, Shanghai 200092, China}
\email{pingli@tongji.edu.cn}
\thanks{Research supported by Program for Young Excellent
Talents in Tongji University.}


\subjclass[2000]{Primary 57R20; Secondary 58J20.}


\keywords{Chern number, index, Dirac operator, signature operator}

\begin{abstract}
Libgober and Wood proved that the Chern number $c_{1}c_{n-1}$ of a
$n$-dimensional compact complex manifold can be determined by its
Hirzebruch $\chi_{y}$-genus. Inspired by the idea of their proof, we
show that, for compact, spin, almost-complex manifolds, more Chern
numbers can be determined by the indices of some twisted Dirac and
signature operators. As a byproduct, we get a divisibility result of
certain characteristic number for such manifolds. Using our method,
we give a direct proof of Libgober-Wood's result, which was
originally proved by induction.
\end{abstract}

\maketitle

\section{Introduction and main results}
 Suppose $(M,J)$ is a compact, almost-complex
$2n$-manifold with a given almost complex structure $J$. This $J$
makes the tangent bundle of $M$ into a $n$-dimensional complex
vector bundle $T_{M}$. Let $c_{i}(M,J)\in H^{2i}(M;\mathbb{Z})$ be
the $i$-th Chern class of $T_{M}$. Suppose we have a formal
factorization of the total Chern class as
follows:$$1+c_{1}(M,J)+\cdots+c_{n}(M,J)=\prod_{i=1}^{n}(1+x_{i}),$$
i.e., $x_{1},\cdots,x_{n}$ are formal Chern roots of $T_{M}$. The
Hirzebruch $\chi_{y}$-genus of $(M,J)$, $\chi_{y}(M,J)$, is defined
by
$$\chi_{y}(M,J)=[\prod_{i=1}^{n}\frac{x_{i}(1+ye^{-x_{i}})}{1-e^{-x_{i}}}][M].$$
Here $[M]$ is the fundamental class of the orientation of $M$
induced by $J$ and $y$ is an indeterminate. If $J$ is specified, we
simply denote $\chi_{y}(M,J)$ by $\chi_{y}(M).$

When the almost complex structure $J$ is integrable, i.e., $M$ is a
(complex) $n$-dimensional compact complex manifold, $\chi_{y}(M)$
can be obtained by
$$\chi^{p}(M)=\sum_{q=0}^{n}(-1)^{q}h^{p,q}(M),~~\chi_{y}(M)=\sum_{p=0}^{n}\chi^{p}(M)\cdot
y^{p},$$ where $h^{p,q}(\cdot)$ is the Hodge number of type $(p,q)$.
This is given by the Hirzebruch-Riemann-Roch Theorem. It was first
proved by Hirzebruch (\cite{Hi}) for projective manifolds, and in
the general case by Atiyah and Singer (\cite{AS}).

The formula \be\label{HRR}\sum_{p=0}^{n}\chi^{p}(M)\cdot
y^{p}=[\prod_{i=1}^{n}\frac{x_{i}(1+ye^{-x_{i}})}{1-e^{-x_{i}}}][M]\ee
implies that $\chi^{p}(M)$ (the index of the Dolbeault complex) can
be expressed as a rationally linear combination of some Chern
numbers of $M$. Conversely, we can address the following question.

\begin{question} For a $n$-dimensional compact complex manifold $M$, given a partition
$\lambda=\lambda_{1}\lambda_{2}\cdots\lambda_{l}$ of weight $n$,
whether the corresponding Chern number
$c_{\lambda_{1}}c_{\lambda_{2}}\cdots c_{\lambda_{l}}[M]$ can be
determined by $\chi^{p}(M)$, or more generally by the indices of
some other elliptic differential operators?
\end{question}

For the simplest case $c_{n}[M]$, the answer is affirmative and
well-known (\cite{Hi}, Theorem
15.8.1):
$$c_{n}[M]=\chi_{y}(M)|_{y=-1}=\sum_{p=0}^{n}(-1)^{p}\chi^{p}(M).$$
The next-to-simplest case is the Chern number $c_{1}c_{n-1}[M]$. The
answer for this Chern number is also affirmative, which was proved
by Libgober and Wood (\cite{LW}, p.141-143):
\be\label{LW}\sum_{p=2}^{n}(-1)^{p}{p\choose
2}\chi^{p}(M)=\frac{n(3n-5)}{24}c_{n}[M]+\frac{1}{12}c_{1}c_{n-1}[M].\ee
The idea of their proof is quite enlightening: expanding both sides
of (\ref{HRR}) at $y=-1$ and comparing the coefficients of the term
$(y+1)^{2}$, then they got (\ref{LW}).

Inspired by this idea, in this paper we will consider the twisted
Dirac operators and signature operators on compact, \emph{spin,
almost-complex} manifolds and show that the Chern numbers $c_{n}$,
$c_{1}c_{n-1}$, $c^{2}_{1}c_{n-2}$ and $c_{2}c_{n-2}$ of such
manifolds can be determined by the indices of these operators.
\begin{remark}
It is worth pointing out that (\ref{LW}) was also obtained later in
(\cite{Sa}, p.144) by Salamon and he majorly applied this result to
hyper-K\"{a}hler manifolds.
\end{remark}

Before stating our main result, we need to fix some definitions and
symbols. Let $M$ be a compact, \emph{almost-complex} $2n$-manifold.
We still use $x_{1},\cdots,x_{n}$ to denote the corresponding formal
Chern roots of the $n$-dimensional complex vector bundle $T_{M}$.
Suppose $E$ is a complex vector bundle over $M$. Set
$$\hat{A}(M, E):=[\textrm{ch}(E)\cdot\prod_{i=1}^{n}\frac{x_{i}/2}{\textrm{sinh}(x_{i}/2)}][M],$$
and
$$L(M, E):=[\textrm{ch}(E)\cdot\prod_{i=1}^{n}\frac{x_{i}(1+e^{-x_{i}})}{1-e^{-x_{i}}}][M],$$
where $\textrm{ch}(E)$ is the Chern character of $E$. The celebrated
Atiyah-Singer index theorem (cf. \cite{HBJ}, p.74 and p.81) states
that $L(M, E)$ is the index of the signature operator twisted by $E$
and when $M$ is \emph{spin}, $\hat{A}(M, E)$ is the index of the
Dirac operator twisted by $E$.

\begin{definition}$$A_{y}(M):=\sum_{p=0}^{n}\hat{A}(M,
\Lambda^{p}(T^{\ast}_{M}))\cdot
y^{p},$$
$$L_{y}(M):=\sum_{p=0}^{n}L(M,
\Lambda^{p}(T^{\ast}_{M}))\cdot y^{p},$$ where
$\Lambda^{p}(T^{\ast}_{M})$ denotes the $p$-th exterior power of the
dual of $T_{M}$.
\end{definition}
Our main result is the following
\begin{theorem}\label{MR}If $M$ is a compact, almost-complex
manifold, then we have
\begin{enumerate}
\item $$\sum_{p=0}^{n}(-1)^{p}\hat{A}(M,
\Lambda^{p}(T^{\ast}_{M}))=c_{n}[M],$$
$$\sum_{p=1}^{n}(-1)^{p}\cdot p\cdot\hat{A}(M,
\Lambda^{p}(T^{\ast}_{M}))=\frac{1}{2}\{nc_{n}[M]+c_{1}c_{n-1}[M]\},$$
$$\sum_{p=2}^{n}(-1)^{p}{p\choose 2}\hat{A}(M,
\Lambda^{p}(T^{\ast}_{M}))=\{\frac{n(3n-5)}{24}c_{n}+\frac{3n-2}{12}c_{1}c_{n-1}+\frac{1}{8}c_{1}^{2}c_{n-2}\}[M];$$
\item $$\sum_{p=0}^{n}(-1)^{p}L(M,
\Lambda^{p}(T^{\ast}_{M}))=2^{n}c_{n}[M],$$
$$\sum_{p=1}^{n}(-1)^{p}\cdot p\cdot L(M,
\Lambda^{p}(T^{\ast}_{M}))=2^{n-1}\{nc_{n}[M]+c_{1}c_{n-1}[M]\},$$
$$\sum_{p=2}^{n}(-1)^{p}{p\choose 2}L(M,
\Lambda^{p}(T^{\ast}_{M}))=2^{n-2}\{\frac{n(3n-5)}{6}c_{n}+\frac{3n-2}{3}c_{1}c_{n-1}+c_{1}^{2}c_{n-2}-c_{2}c_{n-2}\}[M].$$
\end{enumerate}
\end{theorem}
\begin{corollary}
\begin{enumerate}
\item
The Chern numbers $c_{n}[M]$, $c_{1}c_{n-1}[M]$ and
$c^{2}_{1}c_{n-2}[M]$ can be determined by $A_{y}(M)$.

\item
The characteristic numbers $c_{n}[M]$, $c_{1}c_{n-1}[M]$ and
$c_{1}^{2}c_{n-2}[M]-c_{2}c_{n-2}[M]$ can be determined by
$L_{y}(M)$.

\item
When $M$ is \emph{spin}, the Chern numbers $c_{n}[M]$,
$c_{1}c_{n-1}[M]$, $c^{2}_{1}c_{n-2}[M]$ and $c_{2}c_{n-2}[M]$ can
be expressed by using linear combinations of the indices of some
twisted Dirac and signature operators.
\end{enumerate}
\end{corollary}
As remarked in page 143 of \cite{LW}, it was shown by Milnor (cf.
\cite{Hi2}) that every complex cobordism class contains a
non-singular algebraic variety. Milnor (\cite{Mi}) also showed that
two almost-complex manifolds are complex cobordant if and only if
they have the same Chern numbers. Hence Libgober and Wood's result
implies that the characteristic number
$$\frac{n(3n-5)}{24}c_{n}[N]+\frac{1}{12}c_{1}c_{n-1}[N]$$
is always an integer for any compact, \emph{almost-complex}
$2n$-manifold $N$.

Note that the right-hand side of the third equality in Theorem
\ref{MR} is
$$\{\frac{n(3n-5)}{24}c_{n}[M]+\frac{1}{12}c_{1}c_{n-1}[M]\}+\frac{1}{8}\{2(n-1)c_{1}c_{n-1}[M]+c_{1}^{2}c_{n-2}[M]\}.$$
Hence we get the following divisibility result.
\begin{corollary}
For a compact, \emph{spin}, almost-complex manifold $M$, the integer
$2(n-1)c_{1}c_{n-1}[M]+c_{1}^{2}c_{n-2}[M]$ is divisible by 8.
\end{corollary}
\begin{example}
The total Chern class of the complex projective space
$\mathbb{C}P^{n}$ is given by $c(\mathbb{C}P^{n})=(1+g)^{n+1}$,
where $g$ is the standard generator of
$H^{2}(\mathbb{C}P^{n};\mathbb{Z})\cong\mathbb{Z}.$
$\mathbb{C}P^{n}$ is spin if and only if $n$ is odd as
$c_{1}(\mathbb{C}P^{n})=(n+1)g$. Suppose $n=2k+1$. Then
$$2(n-1)c_{1}c_{n-1}[\mathbb{C}P^{n}]+c_{1}^{2}c_{n-2}[\mathbb{C}P^{n}]=8(k+1)^{2}[k(2k+1)+\frac{1}{3}k(k+1)(2k+1)].$$
While it is easy to check that $\mathbb{C}P^{4}$ does not satisfy
this divisibility result.
\end{example}
\section{Proof of the main result}
\begin{lemma}
$$A_{y}(M)=\{\prod_{i=1}^{n}[\frac{x_{i}(1+ye^{-x_{i}(1+y)})}{1-e^{-x_{i}(1+y)}}\cdot e^{-\frac{x_{i}(1+y)}{2}}]\}[M],$$
$$L_{y}(M)=\{\prod_{i=1}^{n}[\frac{x_{i}(1+ye^{-x_{i}(1+y)})}{1-e^{-x_{i}(1+y)}}\cdot (1+e^{-x_{i}(1+y)})]\}[M].$$
\end{lemma}
\begin{proof}
From $c(T_{M})=\prod_{i=1}^{n}(1+x_{i})$ we have (see, for example,
\cite{HBJ}, p.11)
$$c(\Lambda^{p}(T_{M}^{\ast}))=\prod_{1\leq i_{1}<\cdots<i_{p}\leq n}[1-(x_{i_{1}}+\cdots+x_{i_{p}})].$$
Hence
$$ch(\Lambda^{p}(T_{M}^{\ast}))y^{p}=\sum_{1\leq i_{1}<\cdots<i_{p}\leq n}e^{-(x_{i_{1}}+\cdots+x_{i_{p}})}y^{p}=\sum_{1\leq i_{1}<\cdots<i_{p}\leq n}(\prod_{j=1}^{p}ye^{-x_{i_{j}}}).$$
Therefore we have
$$\sum_{p=0}^{n}ch(\Lambda^{p}(T_{M}^{\ast}))y^{p}=\sum_{p=0}^{n}[\sum_{1\leq i_{1}<\cdots<i_{p}\leq n}(\prod_{j=1}^{p}ye^{-x_{i_{j}}})]=\prod_{i=1}^{n}(1+ye^{-x_{i}}).$$
So \be\label{A(M)}
\begin{split}
A_{y}(M)& =\sum_{p=0}^{n}\hat{A}(M, \Lambda^{p}(T^{\ast}_{M}))\cdot
y^{p}\\
&=\{[\sum_{p=0}^{n}ch(\Lambda^{p}(T_{M}^{\ast}))y^{p}]\cdot\prod_{i=1}^{n}\frac{x_{i}/2}{\textrm{sinh}(x_{i}/2)}\}[M]\\
&=\{\prod_{i=1}^{n}[(1+ye^{-x_{i}})\cdot\frac{x_{i}/2}{\textrm{sinh}(x_{i}/2)}]\}[M]\\
&=\{\prod_{i=1}^{n}[\frac{x_{i}(1+ye^{-x_{i}})}{1-e^{-x_{i}}}\cdot
e^{-\frac{x_{i}}{2}}]\}[M].
\end{split}\ee
Since for the evaluation only the homogeneous component of degree
$n$ in the $x_{i}$ enters, then we obtain an additional factor
$(1+y)^{n}$ if we replace $x_{i}$ by $x_{i}(1+y)$ in (\ref{A(M)}).
We therefore obtain:

\be\begin{split}
A_{y}(M)&=\{\frac{1}{(1+y)^{n}}\prod_{i=1}^{n}[\frac{x_{i}(1+y)(1+ye^{-x_{i}(1+y)})}{1-e^{-x_{i}(1+y)}}\cdot
e^{-\frac{x_{i}(1+y)}{2}}]\}[M]\\
&=\{\prod_{i=1}^{n}[\frac{x_{i}(1+ye^{-x_{i}(1+y)})}{1-e^{-x_{i}(1+y)}}\cdot
e^{-\frac{x_{i}(1+y)}{2}}]\}[M].
\end{split}\nonumber\ee
Similarly, \be\begin{split}
L_{y}(M)&=\{\prod_{i=1}^{n}[(1+ye^{-x_{i}})\cdot\frac{x_{i}(1+e^{-x_{i}})}{1-e^{-x_{i}}}]\}[M]\\
&=\{\frac{1}{(1+y)^{n}}\prod_{i=1}^{n}[\frac{x_{i}(1+y)(1+ye^{-x_{i}(1+y)})}{1-e^{-x_{i}(1+y)}}\cdot
(1+e^{-x_{i}(1+y)})]\}[M]\\
&=\{\prod_{i=1}^{n}[\frac{x_{i}(1+ye^{-x_{i}(1+y)})}{1-e^{-x_{i}(1+y)}}\cdot
(1+e^{x_{i}(1+y)})]\}[M].
\end{split}\nonumber\ee
\end{proof}

\begin{lemma}\label{EXP}
Set $z=1+y$. We have
$$A_{y}(M)=\{\prod_{i=1}^{n}[(1+x_{i})-(x_{i}+\frac{1}{2}x_{i}^{2})z+(\frac{11}{24}x_{i}^{2}+\frac{1}{8}x_{i}^{3})z^{2}+\cdots]\}[M],$$
$$L_{y}(M)=\{\prod_{i=1}^{n}[2(1+x_{i})-(2x_{i}+x_{i}^{2})z+(\frac{7}{6}x_{i}^{2}+\frac{1}{2}x_{i}^{3})z^{2}+\cdots]\}[M].$$
\end{lemma}
\begin{proof}
\be
\begin{split}
\frac{x_{i}(1+ye^{-x_{i}(1+y)})}{1-e^{-x_{i}(1+y)}}&=-x_{i}y+\frac{x_{i}(1+y)}{1-e^{-x_{i}(1+y)}}=-x_{i}(z-1)+\frac{x_{i}z}{1-e^{-x_{i}z}}\\
&=-x_{i}(z-1)+(1+\frac{1}{2}x_{i}z+\frac{1}{12}x_{i}^{2}z^{2}+\cdots)\\
&=(1+x_{i})-\frac{1}{2}x_{i}z+\frac{1}{12}x_{i}^{2}z^{2}+\cdots.
\end{split}\nonumber\ee
So we have \be
\begin{split}
A_{y}(M)&=\{\prod_{i=1}^{n}[\frac{x_{i}(1+ye^{-x_{i}(1+y)})}{1-e^{-x_{i}(1+y)}}\cdot
e^{-\frac{x_{i}(1+y)}{2}}]\}[M]\\
&=\{\prod_{i=1}^{n}[(1+x_{i})-\frac{1}{2}x_{i}z+\frac{1}{12}x_{i}^{2}z^{2}+\cdots][1-\frac{1}{2}x_{i}z+\frac{1}{8}x_{i}^{2}z^{2}+\cdots]\}[M]\\
&=\{\prod_{i=1}^{n}[(1+x_{i})-(x_{i}+\frac{1}{2}x_{i}^{2})z+(\frac{11}{24}x_{i}^{2}+\frac{1}{8}x_{i}^{3})z^{2}+\cdots]\}[M].
\end{split}\nonumber\ee
Similarly, \be
\begin{split}
L_{y}(M)&=\{\prod_{i=1}^{n}[\frac{x_{i}(1+ye^{-x_{i}(1+y)})}{1-e^{-x_{i}(1+y)}}\cdot (1+e^{-x_{i}(1+y)})]\}[M]\\
&=\{\prod_{i=1}^{n}[(1+x_{i})-\frac{1}{2}x_{i}z+\frac{1}{12}x_{i}^{2}z^{2}+\cdots][2-x_{i}z+\frac{1}{2}x_{i}^{2}z^{2}+\cdots]\}[M]\\
&=\{\prod_{i=1}^{n}[2(1+x_{i})-(2x_{i}+x_{i}^{2})z+(\frac{7}{6}x_{i}^{2}+\frac{1}{2}x_{i}^{3})z^{2}+\cdots]\}[M].
\end{split}\nonumber\ee
\end{proof}
Let $f(x_{1},\cdots,x_{n})$ be a symmetric polynomial in
$x_{1},\cdots,x_{n}$. Then $f(x_{1},\cdots,x_{n})$ can be expressed
in terms of $c_{1},\cdots,c_{n}$ in a unique way. We use
$h(f(x_{1},\cdots,x_{n}))$ to denote the homogeneous component of
degree $n$ in $f(x_{1},\cdots,x_{n})$. For instance, when
$n=3$,
\be
\begin{split}& h(x_{1}+x_{2}+x_{3}+x_{1}^{2}x_{2}+x_{1}^{2}x_{3}+x_{2}^{2}x_{1}+x_{2}^{2}x_{3}+x_{3}^{2}x_{1}+x_{3}^{2}x_{2})\\
&=x_{1}^{2}x_{2}+x_{1}^{2}x_{3}+x_{2}^{2}x_{1}+x_{2}^{2}x_{3}+x_{3}^{2}x_{1}+x_{3}^{2}x_{2}\\
&=(x_{1}+x_{2}+x_{3})(x_{1}x_{2}+x_{1}x_{3}+x_{2}x_{3})-3x_{1}x_{2}x_{3}\\
&=c_{1}c_{2}-3c_{3}. \end{split}\nonumber\ee

The following lemma is a crucial technical ingredient in the proof
of our main result.
\begin{lemma}\label{maintechnicallemma}~
\begin{enumerate}
\item
$h_{1}:=h\{\sum_{i=1}^{n}[x_{i}\prod_{j\neq i}(1+x_{j})]\}=nc_{n},$
\item
$h_{11}:=h\{\sum_{1\leq i<j\leq n}[x_{i}x_{j}\prod_{k\neq
i,j}(1+x_{k})]\}=\frac{n(n-1)}{2}c_{n},$
\item
$h_{2}:=h\{\sum_{i=1}^{n}[x_{i}^{2}\prod_{j\neq
i}(1+x_{j})]\}=-nc_{n}+c_{1}c_{n-1},$
\item
$h_{12}:=h\{\sum_{1\leq i<j\leq
n}[(x_{i}^{2}x_{j}+x_{i}x_{j}^{2})\prod_{k\neq
i,j}(1+x_{k})]\}=(n-2)(-nc_{n}+c_{1}c_{n-1}),$

\item
$h_{22}:=h\{\sum_{1\leq i<j\leq n}[x_{i}^{2}x_{j}^{2}\prod_{k\neq
i,j}(1+x_{k})]\}=\frac{n(n-3)}{2}c_{n}-(n-2)c_{1}c_{n-1}+c_{2}c_{n-2},$

\item
$h_{3}:=h\{\sum_{i=1}^{n}[x_{i}^{3}\prod_{j\neq
i}(1+x_{j})]\}=nc_{n}-c_{1}c_{n-1}+c_{1}^{2}c_{n-2}-2c_{2}c_{n-2}.$
\end{enumerate}
\end{lemma}

Now we can complete the proof of Theorem \ref{MR} and postpone the
proof of Lemma \ref{maintechnicallemma} to the end of this section.

\begin{proof}
From Lemma \ref{EXP}, the constant coefficient of $A_{y}(M)$ is
$[\prod_{i=1}^{n}(1+x_{i})][M]=c_{n}[M].$

The coefficient of $z$ is
$$\{\sum_{i=1}^{n}[-(x_{i}+\frac{1}{2}x_{i}^{2})\prod_{j\neq i}(1+x_{j})]\}[M]=(-h_{1}-\frac{1}{2}h_{2})[M]=-\frac{1}{2}\{nc_{n}[M]+c_{1}c_{n-1}[M]\}.$$
The coefficient of $z^{2}$ is \be
\begin{split}
&
\{\sum_{i=1}^{n}[(\frac{11}{24}x_{i}^{2}+\frac{1}{8}x_{i}^{3})\prod_{j\neq
i}(1+x_{j})]+\sum_{1\leq i<j\leq
n}[(x_{i}+\frac{1}{2}x_{i}^{2})(x_{j}+\frac{1}{2}x_{j}^{2})\prod_{k\neq
i,j}(1+x_{k})]\}[M]\\
&=(\frac{11}{24}h_{2}+\frac{1}{8}h_{3}+h_{11}+\frac{1}{2}h_{12}+\frac{1}{4}h_{22})[M]\\
&=\{\frac{n(3n-5)}{24}c_{n}+\frac{3n-2}{12}c_{1}c_{n-1}+\frac{1}{8}c_{1}^{2}c_{n-2}\}[M].
\end{split}\nonumber\ee
Similarly, for $L_{y}(M)$, the constant coefficient is
$[2^{n}\prod_{i=1}^{n}(1+x_{i})][M]=2^{n}c_{n}[M].$

The coefficient of $z$ is
$$\{\sum_{i=1}^{n}[-(2x_{i}+x_{i}^{2})\prod_{j\neq i}2(1+x_{j})]\}[M]=(-2^{n}h_{1}-2^{n-1}h_{2})[M]=-2^{n-1}\{nc_{n}[M]+c_{1}c_{n-1}[M]\}.$$

The coefficient of $z^{2}$ is \be
\begin{split}
&
\{\sum_{i=1}^{n}[(\frac{7}{6}x_{i}^{2}+\frac{1}{2}x_{i}^{3})\prod_{j\neq
i}2(1+x_{j})]+\sum_{1\leq i<j\leq
n}[(2x_{i}+x_{i}^{2})(2x_{j}+x_{j}^{2})\prod_{k\neq
i,j}2(1+x_{k})]\}[M]\\
&=(\frac{7\cdot2^{n-2}}{3}h_{2}+2^{n-2}h_{3}+2^{n}h_{11}+2^{n-1}h_{12}+2^{n-2}h_{22})[M]\\
&=2^{n-2}\{\frac{n(3n-5)}{6}c_{n}+\frac{3n-2}{3}c_{1}c_{n-1}+c_{1}^{2}c_{n-2}-c_{2}c_{n-2}\}[M].
\end{split}\nonumber\ee
\end{proof}
It suffices to show Lemma \ref{maintechnicallemma}.

\begin{proof}
(1) and (2) are quite obvious. In the following proof, $\hat{x_{i}}$
means deleting $x_{i}$.

For (3), \be\begin{split}
h_{2}&=\sum_{i=1}^{n}\{h[x_{i}^{2}\prod_{j\neq i}(1+x_{j})]\}\\
&= \sum_{i=1}^{n}(x_{i}\sum_{j\neq i}x_{1}\cdots\hat{x_{j}}\cdots
x_{n})\\
&=\sum_{i=1}^{n}(x_{i}c_{n-1}-c_{n})\\
&=-nc_{n}+c_{1}c_{n-1}.\end{split}\nonumber\ee

For (4),\be\begin{split} h_{12}&=\sum_{1\leq i<j\leq
n}\{h[(x_{i}^{2}x_{j}+x_{i}x_{j}^{2})\prod_{k\neq
i,j}(1+x_{k})]\}\\
&=\sum_{1\leq i<j\leq n}[(x_{i}+x_{j})\sum_{k\neq
i,j}x_{1}\cdots\hat{x_{k}}\cdots x_{n}]\\
&=(n-2)\sum_{i=1}^{n}(x_{i}\sum_{k\neq
i}x_{1}\cdots\hat{x_{k}}\cdots x_{n})\\
&=(n-2)h_{2}\\
&=(n-2)(-nc_{n}+c_{1}c_{n-1}).\end{split}\nonumber\ee

For (5),\be\begin{split} c_{2}c_{n-2}&=(\sum_{1\leq i<j\leq
n}x_{i}x_{j})(\sum_{1\leq k<l\leq
n}x_{1}\cdots\hat{x_{k}}\cdots\hat{x_{l}}\cdots
x_{n})\\
&=\sum_{1\leq i<j\leq n}(x_{i}x_{j}\sum_{1\leq k<l\leq
n}x_{1}\cdots\hat{x_{k}}\cdots\hat{x_{l}}\cdots
x_{n})\\
&=\sum_{1\leq i<j\leq n}[x_{1}x_{2}\cdots
x_{n}+(x_{i}^{2}x_{j}+x_{i}x_{j}^{2})\sum_{k\neq
i,j}x_{1}\cdots\hat{x_{k}}\cdots\hat{x_{i}}\cdots\hat{x_{j}}\cdots
x_{n}+\\
& x_{i}^{2}x_{j}^{2}\sum_{1\leq k<l\leq n;k\neq i,j;l\neq
i,j.}x_{1}\cdots\hat{x_{k}}\cdots\hat{x_{l}}\cdots\hat{x_{i}}\cdots\hat{x_{j}}\cdots
x_{n}]\\
&=\frac{n(n-1)}{2}c_{n}+h_{12}+h_{22}.\end{split}\nonumber\ee
Therefore,
$$h_{22}=c_{2}c_{n-2}-\frac{n(n-1)}{2}c_{n}-h_{12}=\frac{n(n-3)}{2}c_{n}-(n-2)c_{1}c_{n-1}+c_{2}c_{n-2}.$$
For (6),\be\begin{split}
(c_{1}^{2}-2c_{2})c_{n-2}&=(\sum_{i=1}^{n}x_{i}^{2})(\sum_{1\leq
j<k\leq n}x_{1}\cdots\hat{x_{j}}\cdots\hat{x_{k}}\cdots
x_{n})\\
&=\sum_{i=1}^{n}(x_{i}^{2}\sum_{1\leq j<k\leq
n}x_{1}\cdots\hat{x_{j}}\cdots\hat{x_{k}}\cdots
x_{n})\\
&=\sum_{i=1}^{n}[(x_{i}^{3}\sum_{1\leq j<k\leq n;j\neq i;k\neq
i.}x_{1}\cdots\hat{x_{j}}\cdots\hat{x_{k}}\cdots\hat{x_{i}}\cdots
x_{n})+(x_{i}^{2}\sum_{k\neq
i}x_{1}\cdots\hat{x_{k}}\cdots\hat{x_{i}}\cdots
x_{n})]\\
&=h_{3}+h_{2}.\end{split}\nonumber\ee
 Hence
$$h_{3}=(c_{1}^{2}-2c_{2})c_{n-2}-h_{2}=nc_{n}-c_{1}c_{n-1}+c_{1}^{2}c_{n-2}-2c_{2}c_{n-2}.$$
\end{proof}

\section{Concluding remarks}
Libgober and Wood's proof of (\ref{LW}) is by induction (\cite{LW},
p.142, Lemma 2.2). Here, using our method, we can give a quite
direct proof. We have shown that \be\begin{split}
\chi_{y}(M)&=[\prod_{i=1}^{n}\frac{x_{i}(1+ye^{-x_{i}(1+y)})}{1-e^{-x_{i}(1+y)}}][M]\\
&=\{\prod_{i=1}^{n}[(1+x_{i})-\frac{1}{2}x_{i}z+\frac{1}{12}x_{i}^{2}z^{2}+\cdots]\}[M]
\end{split}\nonumber\ee
The coefficient of $z^{2}$ is \be\begin{split}
&\{\sum_{i=1}^{n}[\frac{1}{12}x_{i}^{2}\prod_{j\neq
i}(1+x_{j})]+\sum_{1\leq i<j\leq
n}[\frac{1}{4}x_{i}x_{j}\prod_{k\neq
i,j}(1+x_{k})]\}[M]\\
&=(\frac{1}{12}h_{2}+\frac{1}{4}h_{11})[M]\\
&=\frac{n(3n-5)}{24}c_{n}[M]+\frac{1}{12}c_{1}c_{n-1}[M].
\end{split}\nonumber\ee

It is natural to ask what the coefficients are for higher order
terms $(y+1)^{p}$ $(p\geq 3)$. Unfortunately the coefficients become
very complicated for such terms. Libgober and Wood have given a
detailed remark for higher order terms' coefficients of
$\chi_{y}(M)$ (\cite{LW}, p.144-145). Note that the expression of
$A_{y}(M)$ (resp. $L_{y}(M)$) has an additional factor
$e^{-\frac{x_{i}(1+y)}{2}}$ (resp. $1+e^{x_{i}(1+y)}$) than that of
$\chi_{y}(M)$. Hence we cannot expect that there are \emph{explicit}
expressions of higher order terms' coefficients like Theorem
\ref{MR}.\\

{\bf Acknowledgments.}~I would like to express my sincere gratitude
to my advisors Boju Jiang and Haibao Duan for their consistent
directions and encouragements during the past several years. I also
thank the referees for their careful readings of the earlier version
of this paper. Their comments enhance the quality of the paper.

\bibliographystyle{amsplain}

\end{document}